\documentstyle[leqno,diagram]{article}

\input amssym.def
\input amssym.tex
\font\teneurm=eurm10
\font\seveneurm=eurm7
\font\fiveeurm=eurm5
\newfam\eurmfam
\textfont\eurmfam=\teneurm
\scriptfont\eurmfam=\seveneurm
\scriptscriptfont\eurmfam=\fiveeurm
\def\eurm#1{{\fam\eurmfam\relax#1}}
\newcommand{\bdl}[1]{{\eurm #1}}

\newtheorem{thm}{Theorem}
\newtheorem{prop}[thm]{Proposition}
\newtheorem{lem}[thm]{Lemma}
\newtheorem{eg}[thm]{Example}
\newtheorem{defn}[thm]{Definition}

\newcommand{\End}[1]{\underline{End}( #1)}
\newcommand{\contraction}{\rfloor}

\begin{document}

\title{Spectral covers}
\author{Ron Donagi\thanks{Partially supported by NSF Grant DMS 95-03249 and NSA
Grant MDA 904-92-H3047  }}
\date{}
\maketitle

\section{Introduction} \label{intro}\
\indent Spectral curves arose historically out of the study of differential
equations of Lax type. Following Hitchin's work \cite{H1}, they have acquired a
central role in understanding the moduli spaces of vector bundles and Higgs
bundles on a curve. Simpson's work \cite{S} suggests a similar role for
spectral covers $\widetilde{S}$ of higher dimensional varieties $S$ in moduli
questions for bundles on $S$.

The purpose of these notes is to combine and review  various results about
spectral covers, focusing on the decomposition of their Picards (and the
resulting Prym identities) and the interpretation of a distinguished Prym
component as parameter space for Higgs bundles. Much of this is modeled on
Hitchin's system, which we recall in section \ref{Hitchin}, and on several
other systems based on moduli of Higgs bundles, or vector bundles with twisted
endomorphisms,
 on curves.
By peeling off several layers of data which are not essential for our purpose,
we arrive at the notions of an {\em abstract principal Higgs bundle} and a {\em
cameral}  (roughly, a principal spectral) cover. Following \cite{D3},  this
leads to the statement of the main result, theorem \ref{main}, as an
equivalence between these somewhat abstract  `Higgs' and `spectral' data, valid
over an arbitrary complex variety and for a reductive Lie group $G$. Several
more familiar forms of the equivalence can then be derived in special cases by
adding choices of  representation, value bundle and twisted endomorphism. This
endomorphism is required to be {\em regular}, but not semesimple. Thus the
theory works well even for Higgs bundles which are everywhere nilpotent.  After
touching briefly on the symplectic side of the story In section
\ref{symplectic}, we discuss some of the issues involved in removing the
regularity assumption, as well as some applications and open problems, in sect!
ion \ref{apps}.

This survey is based on talks at  the Vector Bundle Workshop at UCLA (October
92)  and the Orsay meeting (July 92), and earlier talks at Penn, UCLA and MSRI.
 I would like to express my thanks to Rob Lazarsfeld and Arnaud Beauville for
the invitations, and to them and Ching Li Chai, Phillip Griffiths, Nigel
Hitchin, Vasil Kanev, Ludmil Katzarkov, Eyal Markman, Tony Pantev, Emma
Previato and Ed Witten for stimulating and helpful conversations.

We work throughout over $\bf C$ . The total space of a vector bundle  (=locally
free sheaf) $K$ is denoted $\Bbb{K}$. Some more notation:
\begin{tabbing}
bundles of algebras:    \=  \kill
Groups:                          \>$G $           \=$B$          \=$T$
  \=$N $          \=$C$\\
algebras:                        \>$\frak g$   \>$\frak b$  \>$\frak t$
\>$\frak n $   \>$\frak c $\\
Principal bundles:        \>$\cal G$    \>$\cal B$  \>$\cal T$   \>$\cal N $
\>$\cal C$\\
bundles of algebras:    \>\bdl{g}    \>\bdl{b}    \>\bdl{t}    \>\bdl{n}
\>\bdl{c}
\end{tabbing}

\section{Hitchin's system}\
\indent Let  ${\cal M} := {\cal M} _{C}(n,d)$ be the moduli space of stable
vector bundles of rank n and degree d on a smooth projective complex curve $C$.
 It is smooth and quasiprojective of dimension

\begin{equation}
   \tilde{g} := n^{2} (g-1)+1.
\end{equation}

\noindent Its cotangent space at a point  $   E \in {\cal M}  $ is given by

\begin{equation}
  T^*_{E}{\cal M}  :=  H^{0} (  \End{E} \otimes \omega_{C} )
\end{equation}

\noindent where $ \omega_{C}  $ is the canonical bundle of $C$. Our starting
point is:

\begin{thm}[Hitchin\cite{H1}]
       \label{Hitchin}
 The cotangent bundle $T^*{\cal M}$ is an algebraically completely integrable
Hamiltonian system.
\end{thm}

{\em Complete integrability} means that there is a map
\[  h:T^*{\cal M}\longrightarrow B  \]
 to a $\tilde{g}$-dimensional vector space $B$  which is Lagrangian with
respect to the natural symplectic structure on $ T^*{\cal M}$   (i.e. the
tangent spaces to a general fiber $h^{-1}(a)$ over $a \in B$ are maximal
isotropic subspaces with respect to the symplectic form). In this situation one
gets, by contraction with the symplectic form, a trivialization of the tangent
bundle:

\begin{equation}
 T_{h^{-1}(a)}   \stackrel{ \approx}{ \longrightarrow} {\cal O}_{h^{-1}(a)}
\otimes T^{*}_{a} B. \end{equation}

In particular, this produces a family of ({\em `Hamiltonian'}) vector fields on
 $h^{-1}(a)$ which is parametrized by $T^{*}_{a} B$ , and the flows generated
by these on $h^{-1}(a)$ all commute.
{\em Algebraic complete integrability} means additionally that the fibers
$h^{-1}(a)$ are Zariski open subsets of abelian varieties on which the
Hamiltonian flows are linear, i.e. the vector fields are constant.

We describe the idea of the proof in a slightly more general setting, following
\cite{BNR}. Let $K$ be a line bundle on $C$, with total space $ \Bbb K$ . (In
Hitchin's situation, $K$ is $ \omega_{C}  $ and $ \Bbb K$ is $T^*C$ .) A {\em
$K$-valued Higgs bundle} is a pair
\[ (E   \quad,\quad   \phi :E \longrightarrow E \otimes K) \]
consisting of a vector bundle $E$ on $C$ and a $K$-valued endomorphism. One
imposes an appropriate stability condition, and obtains a good moduli space
${\cal M}_K$ parametrizing equivalence classes of  $K$-valued semistable Higgs
bundles, with an open subset ${\cal M}_K^s$ parametrizing isomorphism classes
of stable ones, cf. \cite{S}.

Let $B:=B_K$ be the vector space parametrizing  polynomial maps
\[  p_a :   \Bbb K  \longrightarrow   \Bbb K^{n} \]
of the form
\[ p_{a} (x) = x^{n} + a_{1} x^{n-1} + \cdots + a_{n},   \qquad\qquad   a_{i}
\in H^{0}(K ^{\otimes i}).  \]
in other words,

\begin{equation}
B := \bigoplus_{i=1}^{n} H^{0}(K ^{\otimes i}).
\end{equation}

\noindent The assignment

\begin{equation}
(E, \phi) \longmapsto char(\phi):= \det{(xI-\phi})
\end{equation}

\noindent gives a morphism

\begin{equation}
h_K:{\cal M}_{K}\longrightarrow B_{K}.
\end{equation}

In Hitchin's case,  the desired map $h$ is the restriction of $h_{\omega _{C}
}$ to $T^{*}{\cal M}$, which is an open subset of  ${\cal M}_{\omega_C}^s$.
Note that in this case  $\dim{B}$  is, in Hitchin's words, `somewhat
miraculously' equal to $\tilde{g} = \dim{{\cal M}}$.

The {\em spectral curve}  $\widetilde{C}:= \widetilde{C}_{a}$ defined by $a \in
B_{K}$ is the inverse image in $\Bbb K$ of the $0$-section of  $\Bbb K
^{\otimes n}$ under $  p_a :   \Bbb K  \longrightarrow   \Bbb K^{n} $. It is
finite over $C$ of degree $n$. The general fiber of $h_K$ is given by:

\begin{prop} \cite{BNR}         \label{prop:BNR}
For $a \in B$  with {\em integral} spectral curve $\widetilde{C}_{a}$, there is
a natural equivalence  between isomorphism classes of:
\begin{enumerate}
  \item Rank-$1$, torsion-free sheaves on $\widetilde{C}_{a}$.
  \item  Pairs $(E    \   ,  \    \phi:E \rightarrow E \otimes K)$ with
$char(\phi)=a$.
\end{enumerate}
\end{prop}

When $\widetilde{C}_{a}$ is non-singular, the fiber is thus
$Jac(\widetilde{C}_{a})$,  an abelian variety.
In $T^*\cal M$ the fiber is an open subset of this abelian variety. One checks
that the missing part has codimension $ \geq 2$, so the symplectic form, which
is exact, must  restrict to $0$ on the fibers, completing the proof.

\section  {Some related systems}

\noindent \underline{\bf Polynomial matrices} \nopagebreak

\noindent   One of the earliest appearances of an ACIHS (algebraically
completely integrable Hamiltonian system) was in Jacobi's work on the geodesic
flow on an ellipsoid (or more generally, on a nonsingular quadric in ${\bf
R}^k$). Jacobi discovered that this differential equation, taking place  on the
tangent (=cotangent!)  bundle of the ellipsoid, can be integrated explicitly in
terms of hyperelliptic theta functions. In our language, the total space of the
flow is an ACIHS, fibered by (Zariski open subsets of) hyperelliptic Jacobians.
We are essentially in the special case of
 Proposition \ref{prop:BNR}  where
\[ C={\bf  P}^1, \quad         n=2, \quad        K= {\cal O}_{{\bf P}^1} (k).
\]
A variant of this system appeared in Mumford's solution \cite{Mu}   of the
Schottky problem for hyperelliptic curves.

The extension to all values of n is studied in \cite{B} and, somewhat  more
analytically, in  \cite{AHP} and \cite{AHH}. Beauville considers, for fixed $n$
and $k$, the space $B$ of polynomials:

\begin{equation}
p=y^n + a_1(x)y^{n-1}+\cdots +a_n(x), \quad, \deg{(a_i)}\leq ki
\end{equation}

\noindent
in variables $x$ and $y$.  Each $p$ determines an $n$-sheeted branched cover
$$\widetilde{C}_p \rightarrow {\bf P}^1.$$
The  total space is the space of polynomial matrices

\begin{equation}
M := H^0 ({\bf  P}^1 ,  \End{{\cal O}^{\oplus n}}  \otimes {\cal O}(d) ),
\end{equation}

\noindent
the map $h:M \rightarrow B$ is the characteristic polynomial,
and $M_p$ is the fiber over a given $p \in  B$.
The result is that  for smooth spectral curves $\widetilde{C}_p$,
$ {\bf P}GL(n)$ acts freely and properly on $M_p$;
the quotient is isomorphic to
 $J(\widetilde{C}_p) \smallsetminus \Theta. $
(In order to obtain the entire $J(\widetilde{C}_p) ,$ one must allow all pairs
$(E,\phi)$ with $E$ of given degree, say $0$. Among those, the ones with
 $E\approx  {{\cal O}_{{\bf  P}^1}}^{\oplus n}$
correspond to  the open set
$J(\widetilde{C}_p) \smallsetminus \Theta. $ )
This system is an ACIHS, in a slightly weaker sense than before: instead of a
symplectic structure, it has a {\em Poisson structure}, i.e.  a section $\beta$
of $\wedge^2 T$, such that the $\bf C$-linear sheaf map given by contraction
with $\beta$
$$\begin{array}{ccc}
   {\cal O}   &   \rightarrow    &    {\cal T}   \\
   f                &   \mapsto    &    df  \contraction  \beta
\end{array}$$
sends the Poisson bracket of functions to the bracket of vector fields. Any
Poisson manifold is naturally foliated, with (locally analytic) symplectic
leaves.  For a Poisson ACIHS, we want each leaf to inherit a (symplectic)
ACIHS, so the symplectic foliation should be pulled back via $h$ from a
foliation of the base $B$.

The result of \cite{BNR} suggests that analogous systems should exist when
${\bf P}^1$ is replaced by an arbitrary base curve $C$.  The main point is to
construct the Poisson structure. This was achieved by Bottacin \cite{Bn} and
Markman \cite{M1}, cf. section \ref{symplectic}.
In the case of the polynomial matrices though, everything (the commuting vector
fields, the Poisson structure, etc.) can be written very explicitly. What makes
these explicit results  possible is that  every vector bundle over ${\bf P}^1$
splits. This of course fails in genus $>1$, but for elliptic curves  the moduli
space of vector bundles is still completely understood, so here too the system
can be described explicitly:

  For simplicity, consider vector bundles with fixed determinant.  When the
degree is $0$, the moduli space is a projective space ${\bf P}^{n-1}$ (or more
canonically, the fiber over $0$ of the Abel-Jacobi map
$$C^{[n]} \longrightarrow  J(C) = C.$$
The ACIHS which arises is essentially the Treibich-Verdier theory \cite{TV} of
elliptic solitons. When, on the other hand, the degree is $1$ (or more
generally, relatively prime to $n$), the moduli space is a single point;  the
corresponding system was studied explicitly in \cite{RS}. \\

\noindent \underline{\bf Reductive groups} \nopagebreak

\noindent  In another direction, one can replace the vector bundles by
principal $G$-bundles ${\cal G}$ for any reductive group $G$.  Again, there is
a moduli space ${\cal M}_{G,K}$ parametrizing equivalence classes of semistable
$K$-valued $G$-Higgs bundles, i.e. pairs
$({\cal G}, \phi)$  with $\phi \in K \otimes \bdl{ad}(\cal G)$. The Hitchin map
goes to $$B:=\oplus_{i} H^0(K^{\otimes d_i}),$$ where the $d_i$ are the degrees
of the $f_i$, a basis for the $G$-invariant polynomials on the Lie algebra
$\frak g$:
\[ h: ({\cal G}, \phi) \longrightarrow (f_i (\phi))_{i}.
\]
When $K=\omega_C$, Hitchin showed \cite{H1} that one still gets a completely
integrable system, and that it is algebraically completely integrable for the
classical groups $GL(n), SL(n), SP(n), SO(n).$ The generic fibers are in each
case  (not quite canonically; one must choose various square roots! cf.
sections \ref{reg.ss} and \ref{reg}) isomorphic to abelian varieties given in
terms of the spectral curves $\widetilde{C}$:

\begin{center}
\begin{equation}
\begin{array}{cl}                                      \label{Pryms for groups}
  GL(n)&    \widetilde{C}
                       \mbox{ has degree n over C, the AV is Jac(}
                       \widetilde{C}).    \\
  SL(n)&    \widetilde{C}
                      \mbox{ has degree n over C, the AV is Prym(}
                      \widetilde{C} / C).  \\
  SP(n)&    \widetilde{C}
                      \mbox{ has degree  2n over C and an involution }
                       x  \mapsto -x.  \\
             &    \mbox{ The map factors through the quotient }
                       \overline{C}.   \nonumber \\
             &    \mbox{ The AV is }
                       Prym( \widetilde{C} / \overline{C}).  \nonumber \\
  SO(n)&   \widetilde{C}  \mbox{ has degree  n and an involution , with: }  \\
             &  \bullet \mbox{ a fixed  component,  when n is odd.} \\
             &  \bullet \mbox{ some fixed double points, when n is even.} \\
             &   \mbox{ One must desingularize }
                      \widetilde{C}
                      \mbox{ and the quotient }
                      \overline{C}, \\
             &    \mbox{and ends up with the Prym  of the} \\
             &    \mbox{desingularized double cover.}  \
\end{array}
\end{equation}
\end{center}

The algebraic complete integrability was verified in \cite{KP1} for the
exceptional group $G_2$.
A sketch of the argument for any reductive $G$ is in \cite{BK},  and a complete
proof was given in \cite{F}. We will outline a proof in section
\ref{abelianization} below.   \\

\noindent \underline{\bf Higher dimensions} \nopagebreak

\noindent  Finally, a sweeping extension of the notion of Higgs bundle is
suggested by the work of Simpson  \cite{S}. To him, a Higgs bundle on a
projective variety S is a vector bundle (or principal $G$-bundle \ldots) $E$
with a  {\em symmetric}, $\Omega^1_S$-valued  endomorphism
\[ \phi : E \longrightarrow E \otimes \Omega^1_S.
\]
Here {\em symmetric} means the vanishing of:
\[ \phi\wedge\phi : E \longrightarrow E \otimes \Omega^2_S,
\]
a condition which is obviously vacuous on curves. He constructs a moduli space
for such Higgs bundles (satisfying appropriate stability conditions), and
establishes diffeomorphisms to corresponding moduli spaces of representations
of $\pi_1(S)$ and of connections.

\section {Decomposition of spectral Picards}
\subsection{The question}\
\indent Let $({\cal G},\phi)$ be a $K$-valued principal Higgs bundle on a
complex variety $S$. Each representation
\[   \rho :  G \longrightarrow Aut(V)
\]
determines an associated $K$-valued Higgs bundle
\[   ( {\cal V} := {\cal G} \times^{G} V, \qquad{\rho}(\phi)\ ),
\]
which in turn determines a spectral cover $\widetilde{S}_V \longrightarrow S$.

The question, raised first in \cite{AvM} when $S={\bf P}^1$, is to relate the
Picard varieties of the
$\widetilde{S}_V$ as $V$ varies, and in particular to find pieces common to all
of them. For Adler and van Moerbeke, the motivation was that  many evolution
DEs (of Lax type) can be {\em linearized} on the Jacobians of spectral curves.
This means that the "Liouville tori", which live naturally in the complexified
domain of the DE (and hence are independent of the representation $V$) are
mapped isogenously to their image in $\mbox{Pic}(\widetilde{S}_V)$ for each
nontrivial $V$ ; so one should be able to locate these tori among the pieces
which occur in an isogeny decomposition of each of the
$\mbox{Pic}(\widetilde{S}_V)$. There are many specific examples where a pair of
 abelian varieties constructed from related covers of curves are known to be
isomorphic or isogenous, and some of these lead to important identities among
theta functions.

\begin{eg}
\begin{em}
Take $G=SL(4)$ . The standard representation $V$ gives a branched cover
$\widetilde{S}_V \longrightarrow S$ of degree 4. On the other hand, the
6-dimensional representation $\wedge ^2 V$ (=the standard representation of the
isogenous group $SO(6)$)  gives a cover
$ \stackrel{\approx}{S} \longrightarrow S$  of degree 6, which factors through
an involution:
\[  \stackrel{\approx}{S} \longrightarrow \overline{S} \longrightarrow S.
\]
One has the isogeny decompositions:
\[ Pic \, (\widetilde{S}) \sim Prym(\widetilde{S} / S) \oplus Pic \,(S)
\]
\[ Pic \,(\stackrel{\approx}{S}) \sim
Prym(\stackrel{\approx}{S} / \overline{S})        \oplus
Prym(\overline{S} / S)              \oplus                     Pic \,(S).
\]
It turns out that
\[  Prym(\widetilde{S} / S)            \sim        Prym(\stackrel{\approx}{S} /
\overline{S}) .
\]
For  $S={\bf P}^1$, this is Recillas' {\em trigonal construction} \cite{R}. It
says that every Jacobian of a trigonal curve is the Prym of a double cover of a
tetragonal curve, and vice versa.
\end{em}
\end{eg}

\begin{eg}
\begin{em}
Take $G=SO(8)$ with its standard 8-dimensional representation $V$.  The
spectral cover has degree 8 and factors through an involution,
 $ \stackrel{\approx}{S} \longrightarrow \overline{S} \longrightarrow S.$
The two half-spin representations $V_1, V_2$ yield similar covers
\[ \stackrel{\approx}{S} _i \longrightarrow \overline{S} _i \longrightarrow S,
\qquad i=1,2.
\]
The {\em tetragonal construction} \cite{D1} says that the three Pryms of the
double covers are isomorphic. (These examples, as well as Pantazis' {\em
bigonal construction} and constructions based on some exceptional groups, are
discussed in the context of spectral covers in \cite{K} and \cite{D2}.)
\end{em}
\end{eg}

It turns out in general that there is indeed a distinguished, Prym-like isogeny
component common to all the spectral Picards, on which the solutions to
Lax-type DEs evolve linearly.  This was noticed in some cases already in
\cite{AvM}, and was greatly extended by Kanev's construction of Prym-Tyurin
varieties. (He still needs $S$ to be ${\bf P}^1$ and the spectral cover to have
generic ramification; some of his results apply only to {\em minuscule
representations}.)
Various parts of the general story have been worked out recently by a number of
authors, based on either of two approaches: one, pursued in  \cite{D2,Me,MS},
is to decompose everything according to the action of the Weyl group $W$ and to
look for common pieces; the other, used in \cite{BK,D3,F,Sc}, relies on the
correspondence of spectral data and Higgs bundles . The group-theoretic
approach is described in the rest of this section. We take up the second
method, known as {\em abelianization}, in section~\ref{abelianization}.

\subsection{Decomposition of spectral covers}                    \label{decomp
covers}\
\indent  The decomposition of spectral Picards arises from three sources.
First, the spectral cover for a sum of representations is the union of the
individual covers $\widetilde{S}_V$.  Next, the cover  $\widetilde{S}_V$ for an
irreducible representation is still the union of subcovers
$\widetilde{S}_{\lambda}$  indexed by weight orbits.  And finally, the Picard
of  $\widetilde{S}_{\lambda}$ decomposes into Pryms.
We start with a few observations about the dependence of the covers themselves
on the representation.  The decomposition of the Picards is taken up in the
next subsection. \\

\noindent \underline{\bf Spectral covers} \nopagebreak

\noindent There is an {\em infinite} collection (of  irreducible
representations $V := V_{\mu}$, hence) of  spectral covers $\widetilde{S}_V$,
which can be parametrized by their highest weights $\mu$  in the dominant  Weyl
chamber $\overline{C}$ ,  or equivalently by the $W$-orbit of extremal weights,
in $\Lambda  / W$.  Here $T$ is a maximal torus in $G$, $\Lambda := Hom(T, {\bf
C}^*)$ is the {\em weight lattice } (also called {\em character lattice })  for
$G$, and $W$ is the Weyl group. Each of these $\widetilde{S}_V$ decomposes as
the union of its subcovers $\widetilde{S}_{\lambda}$, parametrizing eigenvalues
in a given $W$-orbit   $W{\lambda}$ . ($\lambda$ runs over the weight-orbits in
$V_{\mu}$.) \\

\noindent \underline{\bf Parabolic covers} \nopagebreak

\noindent There is a {\em finite} collection of covers $\widetilde{S}_P$,
parametrized by the conjugacy classes in $G$ of parabolic subgroups (or
equivalently by arbitrary dimensional faces $F_P$ of the chamber
$\overline{C}$) such that (for general $S$) each  eigenvalue cover
$\widetilde{S}_{\lambda}$ is birational to some  parabolic cover
$\widetilde{S}_{P}$, the one whose open face $F_P$ contains ${\lambda}$.  \\

\noindent \underline{\bf The cameral cover} \nopagebreak

\noindent There is a $W$-Galois cover   $\widetilde{S} \longrightarrow S$ such
that each
 $\widetilde{S}_{P}$ is isomorphic to  $\widetilde{S} / W_P$, where $W_P$ is
the Weyl subgroup of $W$ which stabilizes $F_P$. We call $\widetilde{S}$ the
{\em cameral cover} ,  since, at least generically,  it parametrizes the
chambers determined by $\phi$ (in the duals of the Cartans), or equivalently
the Borel subalgebras containing  $\phi$. This is constructed as follows: There
is a morphism
${\frak g}\longrightarrow {\frak t}/W$ sending $g \in {\frak g}$ to the
conjugacy class of its semisimple part $g_{ss}$. (More precisely, this is
$Spec$ of the composed ring homomorphism
${\bf C} [ {\frak t} ] ^{W}
{ \stackrel{\simeq}{\leftarrow}}
 {\bf C}[{\frak g}]^{G}       \label{t/W}
\hookrightarrow
{\bf C}[{\frak g}]$.)
Taking fiber product with the quotient map ${\frak t}\longrightarrow {\frak
t}/W$, we get the cameral cover ${\tilde{\frak g}}$ of  ${\frak g}$. The
cameral cover $\widetilde{S} \longrightarrow S$ of a $K$-valued principal Higgs
bundle on  $S$  is glued from covers of open subsets in $S$ (on which $K$ and
$\cal G$ are trivialized) which in turn are pullbacks by $\phi$ of
${\tilde{\frak g}}  \longrightarrow  {\frak g} $.

\subsection{Decomposition of spectral Picards}\
\indent The decomposition of  the Picard varieties of spectral covers can be
described as follows:\\

\noindent \underline{\bf The cameral Picard} \nopagebreak

\noindent  From each isomorphism class of irreducible $W$-representations,
choose an integral representative $\Lambda _i$. (This can always be done, for
Weyl groups.)  The group ring
${\bf Z} [W]$  which acts on $Pic(\widetilde{S}) $ has an isogeny
decomposition:

\begin{equation}\label{regular rep}
{\bf Z} [W] \sim \oplus _i \Lambda _i \otimes_{\bf Z} \Lambda _i^{*},
\end{equation}

\noindent
which is just the decomposition of the regular representation. There is a
corresponding isotypic decomposition:

\begin{equation}\label{cameral Pic decomposition}
Pic(\widetilde{S}) \sim \oplus _i \Lambda _i \otimes_{\bf Z} Prym_{\Lambda
_i}(\widetilde{S}),
\end{equation}

\noindent
where

\begin{equation}\label{def of Prym_lambda}
Prym_{\Lambda _i}(\widetilde{S} ):= Hom_W (\Lambda _i , Pic(\widetilde{S})).
\end{equation}\\

\noindent \underline{\bf Parabolic Picards} \nopagebreak

\noindent There are at least three reasonable ways of obtaining an isogeny
decomposition of $Pic(\widetilde{S}_P) $, for a parabolic subgroup $P \subset
G$:
\begin{itemize}

\item The `Hecke' ring $Corr_P$ of correspondences on $\widetilde{S}_P$ over
$S$ acts on  $Pic(\widetilde{S}_P) $, so every irreducible integral
representation $M$ of $Corr_P$ determines a generalized Prym
$$   Hom_{Corr_P} (M, Pic(\widetilde{S}_P)),   $$
and we obtain an isotypic decomposition of $Pic(\widetilde{S}_P)$ as before.

\item $Pic(\widetilde{S}_P)$ maps, with torsion kernel, to
$Pic(\widetilde{S})$, so we obtain a decomposition of the former by
intersecting its image with the isotypic components
$\Lambda _i \otimes_{\bf Z} Prym_{\Lambda _i}(\widetilde{S})$ of the latter.

\item Since $\widetilde{S}_P$ is the cover of $S$ {\em associated} to the
$W$-cover $\widetilde{S}$ via the permutation representation ${\bf Z} [W_P
\backslash W]$ of $W$, we get an isogeny decomposition of
$Pic(\widetilde{S}_P)$ indexed by the irreducible representations in
${\bf Z} [W_P \backslash W]$.

\end{itemize}

It turns out (\cite{D2},section 6) that all three decompositions agree
and can be given explicitly as

\begin{equation}
\label{multiplicity spaces}
\oplus _i M _i \otimes Prym_{\Lambda _i}(\widetilde{S}) \subset
\oplus _i \Lambda _i \otimes Prym_{\Lambda _i}(\widetilde{S}),\qquad
  M_i := (\Lambda_i)^{W_P}.
\end{equation}

\noindent \underline{\bf Spectral Picards} \nopagebreak

\noindent  To obtain the decomposition of the Picards of the original  covers
$\widetilde{S}_V$ or
 $\widetilde{S}_{\lambda}$, we need, in addition to the decomposition of
$Pic(\widetilde{S}_P)$, some information on the singularities. These can arise
from two separate sources:
\begin{description}
\item[Accidental singularities of the $\widetilde{S}_{\lambda}$. ]
For a sufficiently general Higgs bundle, and for a weight $\lambda$ in the
interior of the face $F_P$ of the Weyl chamber  $\overline{C}$, the natural
map:
$$    i_{\lambda}: \widetilde{S}_P\longrightarrow  \widetilde{S}_{\lambda}   $$
is birational. For the {\em standard} representations of the classical groups
of types
 $A_n, B_n$ or $C_n$, this {\em is} an isomorphism. But for general ${\lambda}$
it is {\em not}: In order for $i_{\lambda}$ to be an isomorphism, ${\lambda}$
must be a multiple of a fundamental weight, cf. \cite{D2}, lemma 4.2. In fact,
the list of fundamental weights for which this happens is quite short; for the
classical groups we have only:  $\omega_1$ for  $A_n, B_n$ and $C_n$,
$\omega_n$  (the dual representation) for $A_n$, and $\omega_2$ for $B_2$. Note
that for $D_n$ the list is {\em empty}. In particular, the covers produced by
the standard representation of $SO(2n)$ are singular; this fact, noticed by
Hitchin In \cite{H1},  explains the need for desingularization in his
result~(\ref{Pryms for groups}).

\item[Gluing the $\widetilde{S}_{V}$. ]
In addition to the singularities of each $i_{\lambda}$, there are the
singularities created by the gluing map $\amalg_{\lambda}
\widetilde{S}_{\lambda} \longrightarrow  \widetilde{S}_V$. This makes explicit
formulas somewhat simpler in the case, studied by Kanev \cite{K}, of {\em
minuscule} representations, i.e. representations whose weights form a single
$W$-orbit.  These singularities account, for instance,  for the
desingularization required in the $SO(2n+1)$ case in
(\ref{Pryms for groups}).

\end{description}

\subsection{The distinguished Prym}   \label{distinguished}\
\indent Combining much of the above, the Adler--van Moerbeke problem of finding
a component common to  the $Pic(\widetilde{S}_V)$ for all non-trivial $V$
translates into: \\

\begin{em}
Find the irreducible representations
$\Lambda_i $  of $W$   which occur in  ${\bf Z} [W_P \backslash W] $
with positive multiplicity for all proper Weyl subgroups
$W_P \subsetneqq W.$
\end{em} \\

By Frobenius reciprocity, or (\ref{multiplicity spaces}), this is equivalent to
\\

\begin{em}
Find the irreducible representations
$\Lambda_i $
  of  W such that  for every  proper Weyl subgroup
$W_P \subsetneqq W, $
 the space of invariants
$M_i := (\Lambda_i)^{W_P} $
 is non-zero.
\end{em} \\

One solution is now obvious: the {\em{reflection representation}} of $W$ acting
on the weight lattice $\Lambda$ has this property. In fact,
$\Lambda^{W_P}$ in this case is just the face $F_P$ of $\overline{C}$. The
corresponding component  $Prym_{\Lambda }(\widetilde{S})$ , is called {\em{the
distinguished Prym}.} We will see in section \ref{abelianization} that its
points correspond, modulo some corrections,  to Higgs bundles.

For the classical groups, this turns out to be the only common component. For
$G_2$ and $E_6$ it turns out (\cite{D2}, section 6) that a second common
component exists. The geometric significance of points in these components is
not known. As far as I know, the only component other than the distinguished
Prym which has arisen `in nature' is the one associated to the 1-dimensional
sign representation of $W$, cf. section \ref{apps} and \cite{KP2}.

\section {Abelianization}\label{abelianization}

\subsection{Abstract vs. $K$-valued objects}\
\indent We want to describe the abelianization procedure in a somewhat abstract
setting, as an equivalence between {\em{principal Higgs bundles}} and certain
{\em spectral data}.
Once we fix a {\em{values}} vector bundle $K$, we obtain an equivalence between
{\em $K$-valued principal Higgs bundles} and {\em K-valued spectral data}.
Similarly,
the choice of a representation $V$ of $G$ will determine an equivalence of
{\em $K$-valued Higgs bundles} (of a given representation type) with $K$-valued
spectral data.

As our model of a $W$-cover we take the natural quotient map
$$G/T \longrightarrow G/N
$$
and its partial compactification

\begin{equation}
\overline{G/T} \longrightarrow \overline{G/N}.      \label{partial
compactification}
\end{equation}

Here $T \subset G$ is a maximal torus, and $N$ is its normalizer in $G$.
The quotient $G/N$ parametrizes maximal tori (=Cartan subalgebras) $\frak{t}$
in $\frak{g}$,
while $G/T$ parametrizes pairs ${\frak t \subset \frak b}$
with ${\frak b \subset \frak g}$ a Borel subalgebra.
An element $x \in {\frak g}$ is {\em regular} if the dimension of its
centralizer
${\frak c \subset \frak g}$ equals $\dim{T}$ (=the rank of $\frak{g}$). The
partial compactifications
$ \overline{G/N}$ and $ \overline{G/T}$  parametrize regular centralizers
${\frak c }$ and pairs ${\frak c \subset \frak b}$, respectively.

In constructing the cameral cover in section  \ref{t/W}, we used the $W$-cover
$\frak t \longrightarrow \frak t / W$ and its  pullback cover ${
\widetilde{\frak g} \longrightarrow \frak g}$.
Over the open subset $\frak g_{reg}$ of regular elements, the same cover is
obtained by pulling back (\ref{partial compactification}) via the map
$\alpha : \frak g_{reg}  \longrightarrow \overline{G/N}$ sending an element to
its centralizer:

\begin{equation}
                       \label{commutes}
\begin{array}{lccccc}
\frak t & \longleftarrow & \widetilde{\frak g}_{{reg}} & \longrightarrow &
\overline{G/T}  & \\
\downarrow & &\downarrow & & \downarrow  & \\
\frak t  /W & \longleftarrow & {\frak g}_{{reg}} &
\stackrel{\alpha}{\longrightarrow} & \overline{G/N} &.
\end{array}
\end{equation}

When working with $K$-valued objects, it is usually more convenient to work
with the left hand side of  (\ref{commutes}), i.e. with eigen{\em values}. When
working with the abstract objects, this is unavailable,  so we are forced to
work with the eigen{\em vectors},
or the right hand side of  (\ref{commutes}). Thus:

\begin{defn}
An abstract {\em cameral cover} of $S$ is a finite morphism $\widetilde{S}
\longrightarrow S$
with $W$-action, which locally (etale) in $S$ is a pullback of (\ref{partial
compactification}). \\
\end{defn}

\begin{defn}
A {\em $K$-valued cameral cover}  ($K$ is a vector bundle on $S$) consists of a
cameral cover  $\pi :  \widetilde{S} \longrightarrow S$  together with an
$S$-morphism
\begin{equation}
   \widetilde{S} \times \Lambda  \longrightarrow \Bbb{K}   \label{K-values}
\end{equation}
which is $W$-invariant ($W$ acts on $ \widetilde{S} , \Lambda,$ hence
diagonally on
$\widetilde{S} \times \Lambda $ ) and linear in $\Lambda$. \\
\end{defn}

We note that a cameral cover $\widetilde{S}$ determines quotients
$\widetilde{S}_P$ for parabolic subgroups $P \subset G$. A $K$-valued cameral
cover determines additionally the $\widetilde{S}_{\lambda}$ for $\lambda \in
\Lambda$, as images in $\Bbb{K}$ of
$\widetilde{S} \times \{ \lambda \}$.  The data of (\ref{K-values}) is
equivalent to a $W$-equivariant map $\widetilde{S} \longrightarrow
\frak{t}\otimes_{\bf C} K.$

\begin{defn}                    \label{princHiggs}
A $G$-principal Higgs bundle on $S$ is a pair ($\cal{G}, \bdl{c})$ with
$\cal{G}$ a principal $G$-bundle and  $\bdl{c} \subset ad(\cal{G})$ a subbundle
of regular centralizers.
\\
\end{defn}

\begin{defn}
A $K$-valued $G$-principal Higgs bundle consists of  $( \cal{G}, \bdl{c} )$
as above together with a section $\varphi$ of $\bdl{c} \otimes K$.
\end{defn}

A principal Higgs bundle $(\cal{G}, \bdl{c})$  determines a cameral cover
$\widetilde{S}\longrightarrow S$ and a homomorphism $\Lambda \longrightarrow
\mbox{Pic}(\widetilde{S}).$ Let $F$ be a parameter space for Higgs bundles with
a given $\widetilde{S}$. Each non-zero $\lambda \in \Lambda$ gives a
non-trivial map
$F\longrightarrow \mbox{Pic}(\widetilde{S})$. For $\lambda$ in a face $F_P$ of
$\overline{C}$, this factors through $\mbox{Pic}(\widetilde{S}_P)$. The
discussion in section \ref{distinguished} now suggests that $F$ should be given
roughly by the distinguished Prym,
$$ Hom_W (\Lambda , \mbox{Pic}(\widetilde{S})).
$$
It turns out that this guess needs two corrections. The first correction
involves restricting to a coset of a subgroup; the need for this is visible
even in the simplest case where
$\widetilde{S}$
 is etale over
$S$,
so
$(\cal{G}, \bdl{c})$
is everywhere regular and semisimple
(i.e.
$ \bdl{c}$
 is a bundle of Cartans.)
The second correction involves a twist along the ramification of
$\widetilde{S}$
over
$S$.
We explain these in the next two subsections.

\subsection{The regular semisimple case: the shift}       \label{reg.ss}

\begin{eg}   \label{unramified}
\begin{em}
Fix a smooth projective curve
$C$
and a line bundle
$K \in \mbox{Pic}(C)$
such that
$K^{\otimes 2} \approx  \cal{O}_C.$
This determines an etale double cover
$\pi : \widetilde{C} \longrightarrow C$
with involution
$i$,
and homomorphisms

\begin{center}
$\begin{array}{cccccc}
 \pi^{*}          &:&  \mbox{Pic}(C)                     &\longrightarrow
&\mbox{Pic}(\widetilde{C}) &, \\
 \mbox{Nm} &:&  \mbox{Pic}(\widetilde{C}) &\longrightarrow &\mbox{Pic}(C)
              &, \\
   i^{*}            &:& \mbox{Pic}(\widetilde{C}) &\longrightarrow
&\mbox{Pic}(\widetilde{C}) &,
\end{array}$
\end{center}

satisfying
$$   1+i^{*} = \pi^* \circ \mbox{Nm}.
$$

\begin{itemize}
\item For
$G = GL(2)$
we have
$\Lambda = \bf{Z} \oplus \bf{Z}$,
and
$W = {\cal{S}}_{2}$
permutes the summands, so
$$    Hom_W (\Lambda , \mbox{Pic}(\widetilde{C}))  \approx
\mbox{Pic}(\widetilde{C}).
$$
And indeed, the Higgs bundles corresponding to
$\widetilde{C}$
are parametrized by
$\mbox{Pic}(\widetilde{C})$:
send
$L \in \mbox{Pic}(\widetilde{C})$
to
$(\cal{G}, \bdl{c})$,
where
$\cal{G}$
has associated rank-2 vector bundle
${\cal V} := \pi_* L$,
and
$ \bdl{c} \subset \End{{\cal{V}}}$
is
$\pi_* {\cal O}_{\widetilde{C}}.$
\item On the other hand, for
$G=SL(2)$
we have
$\Lambda=\bf{Z}$
and
$W={\cal{S}}_2$
acts by
$\pm 1$,
so
$$   Hom_W (\Lambda , \mbox{Pic}(\widetilde{S}))  \approx
\{L \in \mbox{Pic}(\widetilde{C})\  | \  i^*L \approx L^{-1} \}
= \mbox{ker}(1+i^*).
$$
This group has 4 connected components. The subgroup
$\mbox{ker(Nm)}$
consists of 2 of these. The connected component of 0 is the classical Prym
variety, cf. \cite{MuPrym}. Now the Higgs bundles correspond, via the above
bijection
$L\mapsto \pi_*L$,
to
$$\{L \in \mbox{Pic}(\widetilde{C}) \ |\ \det (\pi_*L) \approx {\cal O}_C \} =
{\mbox{Nm}}^{-1}(K).
$$
Thus they form the {\em non-zero} coset of the subgroup
$\mbox{ker(Nm)}$.
(If we return to a higher dimensional
$S$, it is possible for $K$ not to be in the image of
$\mbox{Nm}$,
so there may be {\em no}
$SL(2)$-Higgs bundles corresponding to such a cover.)
\end{itemize}
\end{em}
\end{eg}

This example generalizes to all
$G$,
as follows. The equivalence classes of extensions
$$1 \longrightarrow  T \longrightarrow  N'  \longrightarrow  W \longrightarrow
1
$$
(in which the action of $W$ on $T$ is the standard one) are parametrized by the
group cohomology
$H^2(W,T)$.
Here the 0 element corresponds to the semidirect product . The class
$[N] \in H^2(W,T)$
of the normalizer $N$ of $T$ in $G$ may be 0, as it is for
$G=GL(n)  ,  {\bf P}GL(n)  ,  SL(2n+1)  $;
or not, as for
$G=SL(2n)$.

Assume first, for simplicity, that
$S,\widetilde{S}$
are connected and projective. There is then a natural group homomorphism
\begin{equation}
                                   \label{c}
c: Hom_W (\Lambda , \mbox{Pic}(\widetilde{S}))\longrightarrow H^2(W,T).
\end{equation}
Algebraically, this is an edge homomorphism for the Grothendieck spectral
sequence of equivariant cohomology, which gives the exact sequence

\begin{equation}
                 \label{c-edge}\qquad
0 \longrightarrow                            H^1(W,T)
\longrightarrow                                H^1(S,{\cal{C}})
\longrightarrow                                Hom_W (\Lambda ,
\mbox{Pic}(\widetilde{S}))
\stackrel{c}{\longrightarrow}         H^2(W,T).
\end{equation}
where
 ${\cal{C}} := \widetilde{S} \times _W T.$
Geometrically, this expresses a {\em Mumford group} construction: giving
${\cal{L}} \in \mbox{Hom}(\Lambda,\mbox{Pic}(\widetilde{S}))$
is equivalent to giving a principal $T$-bundle
$\cal T$
over
$\widetilde{S}$;
for
${\cal{L}} \in \mbox{Hom}_W(\Lambda,\mbox{Pic}(\widetilde{S}))$,
$c({\cal{L}})$
is the class in
$H^2(W,T)$
of the group
$N'$
of automorphisms of
$\cal T$
which commute with the action on
$\widetilde{S}$
of some
$w \in W$.

To remove the restriction on
$S, \widetilde{S}$,
we need to replace each occurrence of $T$ in (\ref{c}, \ref{c-edge}) by
$\Gamma (\widetilde{S}, T)$,
the global sections of the trivial bundle on
$\widetilde{S}$
with fiber $T$. The natural map
$H^2(W,T) \longrightarrow H^2(W,\Gamma (\widetilde{S}, T))$
allows us to think of
$[N]$
as an element of
$H^2(W,\Gamma (\widetilde{S}, T))$.

\begin{prop} \cite{D3}
Fix an etale $W$-cover
$\pi: \widetilde{S}\longrightarrow S$.
The following data are equivalent:

\begin{enumerate}
\item Principal $G$-Higgs bundles
$(\cal{G}, \bdl{c})$
with cameral cover
$\widetilde{S}$.
\item Principal $N$-bundles
$\cal N$
over $S$ whose quotient by $T$ is
$\widetilde{S}.$
\item $W$-equivariant homomorphisms
${\cal{L}} : \Lambda \longrightarrow \mbox{Pic}(\widetilde{S})$
with
$c({\cal L}) = [N] \in H^2(W,\Gamma (\widetilde{S}, T))$.
\end{enumerate}

\end{prop}

We observe that while the shifted objects correspond to Higgs bundles,
the unshifted objects
$$
{\cal{L}} \in \mbox{Hom}_W(\Lambda,\mbox{Pic}(\widetilde{S})),  \qquad  c({\cal
L})=0
$$

\noindent
come from the $\cal C$-torsers in $H^1(S, {\cal C} ).$

\subsection{The regular case: the twist along the ramification}
\label{reg}

\begin{eg}                               \label{ramified}
\begin{em}
Modify example \ref{unramified} by letting
$K \in \mbox{Pic}(C) $
be arbitrary, and choose a section $b$ of
$K ^{\otimes 2}$
which vanishes on a simple divisor
$B \subset C$.
We get a double cover
$\pi : \widetilde{C} \longrightarrow C$
branched along $B$, ramified along a divisor
$$
R \subset \widetilde{C}, \quad  \pi(R)=B.
$$
Via
$L\mapsto \pi_*L$,
the Higgs bundles still correspond to
$$\{L \in \mbox{Pic}(\widetilde{C}) \ |\ \det (\pi_*L) \approx {\cal O}_C \} =
{\mbox{Nm}}^{-1}(K).
$$
But this is no longer in
$   Hom_W (\Lambda , \mbox{Pic}(\widetilde{S}))$;
rather, the line bundles in question satisfy

\begin{equation}
                   \label{SL(2) twist}
i^*L \approx L^{-1}(R).
\end{equation}

\end{em}
\end{eg}

For arbitrary $G$, let
$\Phi$
denote the root system and
$\Phi^+$
the set of positive roots. There is a decomposition
$$   \overline{G/T} \  \smallsetminus \  G/T   = \bigcup _{\alpha \in
\Phi^+}R_{\alpha}
$$
of the boundary into components, with
$R_{\alpha}$
the fixed locus of the reflection
$\sigma_{\alpha}$
in
$\alpha$.
(Via (\ref{commutes}),  these correspond to the complexified walls in
$\frak t$.)
Thus each cameral cover
$\widetilde{S} \longrightarrow S$
comes with a natural set of (Cartier)  {\em ramification divisors}, which we
still denote
$R_{\alpha}, \quad  \alpha \in \Phi^+.$

For
$w \in W$,
set
$$  F_w := \left\{ \alpha \in \Phi^+ \ | \ w^{-1} \alpha \in \Phi^- \right\}
= \Phi^+ \cap w \Phi^-,
$$
and choose a $W$-invariant form
$\langle , \rangle$
  on
$\Lambda$.
We consider the variety
$$   Hom_{W,R} (\Lambda , \mbox{Pic}(\widetilde{S}))
$$
of $R$-twisted $W$-equivariant homomorphisms, i.e. homomorphisms
$\cal L$
satisfying

\begin{equation}      \qquad
                                \label{G twist}
w^*{\cal L}(\lambda) \approx
 {\cal L}(w\lambda)\left( \sum_{\alpha \in F_w}{
{\langle-2\alpha,w\lambda \rangle \over \langle \alpha ,\alpha \rangle}
R_{\alpha}
} \right) , \qquad  \lambda \in \Lambda, \quad w \in W.
\end{equation}

This turns out to be the correct analogue of  (\ref{SL(2) twist}). (E.g. for a
reflection
$w=\sigma_{\alpha}$,
\quad $F_w$
is
$\left\{ \alpha \right\}$,
so this gives
$ w^*{\cal L}(\lambda) \approx
 {\cal L}(w\lambda)\left(
{{\langle\alpha,2\lambda \rangle \over \langle \alpha,\alpha \rangle}
R_{\alpha}}
 \right),$
which specializes to (\ref{SL(2) twist}).)  As before, there is a class map

\begin{equation}
                                   \label{c,R}
c: Hom_{W,R} (\Lambda , \mbox{Pic}(\widetilde{S}))\longrightarrow
H^2(W,\  \Gamma (\widetilde{S}, T))
\end{equation}

\noindent
which can be described via a Mumford-group construction.

To understand this twist, consider the formal object

\begin{center}
$\begin{array}{cccc}
{1 \over 2} \mbox{Ram}: & \Lambda & \longrightarrow & {\bf Q}\otimes
\mbox{Pic}\widetilde{S},  \\
                               & \lambda & \longmapsto &
\sum_{ ( \alpha \in {\Phi^+} ) }{{\langle\alpha,\lambda \rangle \over \langle
\alpha,\alpha \rangle} R_{\alpha}}.
\end{array}$
\end{center}
In an obvious sense, a principal $T$-bundle
$\cal T$
on
$\widetilde{S}$
(or a homomorphism
${\cal L}: \Lambda \longrightarrow \mbox{Pic}(\widetilde{S})$)
is $R$-twisted $W$-equivariant if and only if
${\cal T} (-{1 \over 2} Ram)$
is $W$-equivariant, i.e. if
${\cal T}$
 and
${1 \over 2} Ram$
transform the same way under $W$.
The problem with this is that
${1 \over 2} Ram$
itself does not make sense as a $T$-bundle, because the coefficients
${\langle\alpha,\lambda\rangle \over \langle\alpha,\alpha\rangle} $
are not integers. (This argument shows that if
$Hom_{W,R} (\Lambda , \mbox{Pic}(\widetilde{S}))$
is non-empty, it is a torser over the untwisted
$Hom_{W} (\Lambda , \mbox{Pic}(\widetilde{S}))$.)

\begin{thm} \cite{D3}
  \label{main}
For a cameral cover
$\widetilde{S} \longrightarrow S$,
the following data are equivalent: \\
(1) $G$-principal Higgs bundles with cameral cover
$\widetilde{S}$. \\
(2) $R$-twisted $W$-equivariant homomorphisms
${\cal L} \in c^{-1}([N]).$
\end{thm}

The theorem has an essentially local nature, as there is no requirement that
$S$ be, say, projective. We also do not need the condition of generic behavior
near the ramification, which appears in \cite{F, Me, Sc}. Thus we may consider
an extreme case, where
$\widetilde{S}$
is `everywhere ramified':

\begin{eg}\begin{em}
In example \ref{ramified}, take the section
$b=0$.
The resulting cover
$\widetilde{C}$
is a `ribbon', or length-2 non-reduced structure on $C$: it is the length-2
neighborhood of $C$ in
$\Bbb{K}$.
The SL(2)-Higgs bundles
$({\cal G},\bdl{c})$
for this
$\widetilde{C}$
have an everywhere nilpotent
$\bdl{c}$,
so the vector bundle
${\cal V} := {\cal G} \times^{SL(2)} V  \approx \pi_* L$
 (where $V$ is the standard 2-dimensional representation) fits in an exact
sequence
$$   0 \longrightarrow  {\cal S} \longrightarrow  {\cal V} \longrightarrow
{\cal Q} \longrightarrow 0
$$
with
${\cal S} \otimes K \approx {\cal Q}.$
Such data are specified by the line bundle
${\cal Q}$,
satisfying
${\cal Q}^{\otimes 2} \approx K$,
and an extension class in
$\mbox{Ext}^1({\cal Q}, {\cal S}) \approx H^1(K^{-1})$.
The kernel of the restriction map
$ \mbox{Pic}(\widetilde{C}) \longrightarrow \mbox{Pic}(C) $
is also given by
$H^1(K^{-1})$
(use the exact sequence
$0 \longrightarrow K^{-1} \longrightarrow  \pi_*{\cal
O}_{\widetilde{C}}^{\times}
   \longrightarrow {\cal O}_C^{\times} \longrightarrow  0$),
and the $R$-twist produces the required square roots of $K$. (For more details
on the nilpotent locus, cf. \cite{L} and \cite{DEL}.)
\end{em}\end{eg}

\subsection{Adding values and representations}\
\indent  Fix a vector bundle $K$, and consider the moduli space $ {\cal
M}_{S,G,K} $ of   $K$-valued $G$-principal Higgs bundles on $S$. (It can be
constructed as in Simpson's \cite{S}, even though the objects we need to
parametrize are slightly different than his. In this subsection we outline a
direct construction.)
It comes with a Hitchin map:

\begin{equation}
  \label{BigHitchin}
  h:   {\cal M}_{S,G,K}  \longrightarrow  B_K
\end{equation}

\noindent where $B := B_K$ parametrizes all possible Hitchin data. Theorem
\ref{main} gives a precise description of the fibers of this map, independent
of the values bundle $K$. This leaves us with the relatively minor task of
describing, for each $K$, the corresponding  base, i.e. the closed subvariety
$B_s$ of $B$ parametrizing {\em split} Hitchin data, or $K$-valued cameral
covers.   The point is that  Higgs bundles satisfy a symmetry condition, which
in Simpson's setup is
$$    \varphi  \wedge \varphi = 0,
$$
and is built into our definition \ref{princHiggs} through the assumption that
\bdl{c} is regular, hence abelian. Since commuting operators have common
eigenvectors, this gives a splitness condition on the Hitchin data, which we
describe below.  (When $K$ is a line bundle, the condition is vacuous, $B_s =
B$.) The upshot is:

\begin{lem}
           \label{parametrization}
The following data are equivalent: \\
(a) A $K$-valued cameral cover of $S$. \\
(b) A split, graded homomorphism
     $R{\bf \dot{\ }} \longrightarrow {Sym}{\bf \dot{\ }}K.$ \\
(c) A split Hitchin datum
$b \in B_s$.
\end{lem}

Here $R{\bf \dot{\ }}$ is the graded ring of $W$-invariant polynomials on
$\frak t$:

\begin{equation}
R{\bf \dot{\ }} := (\mbox{Sym}{\bf \dot{\ }} {\frak t}^*)^W
\approx {\bf C}[\sigma_1,\ldots,\sigma_l],  \qquad  \deg (\sigma_i) = d_i
\end{equation}

\noindent
where
$l := \mbox{Rank}({\frak g})$
and the
$\sigma_i$
form a basis for the $W$-invariant polynomials. The Hitchin base is the vector
space
$$  B :=  B_K  :=  \oplus _{i=1}^l H^0(S, {Sym}^{d_i}K)
\approx \mbox{Hom}(R{\bf \dot{\ }},\mbox{Sym}{\bf \dot{\ }}K).
$$

\noindent For each
$\lambda \in \Lambda$
(or
$\lambda \in {\frak t}^*$,
for that matter), the expression in an indeterminate $x$:

\begin{equation}
               \label{rep poly}
q_{\lambda}(x,t) := \prod_{w \in W}{(x-w\lambda(t))}, \qquad  t \in {\frak t},
\end{equation}
is $W$-invariant (as a function of $t$), so it defines an element
$q_{\lambda}(x) \in R{\bf \dot{\ }}[x].$
A Hitchin datum
$b \in B \approx \mbox{Hom}(R{\bf \dot{\ }},\mbox{Sym}{\bf \dot{\ }}K)$
sends this to
$$  q_{\lambda,b}(x) \in \mbox{Sym}\dot{\ }(K)[x].
$$
We say that $b$ is {\em split} if, at each point of $S$ and for each
$\lambda$,
the polynomial
$q_{\lambda,b}(x)$
factors completely, into terms linear in $x$.

We note that, for $\lambda$ in the interior of $C$ (the positive Weyl chamber),
$q_{\lambda,b}$
gives the equation in
$\Bbb K$
of the spectral cover
$\widetilde{S}_{\lambda}$
of section (\ref{decomp covers}):
$q_{\lambda,b}$ gives a morphism
$\Bbb K  \longrightarrow \mbox{Sym}^N \Bbb K$,
where $N:=\#W$, and $\widetilde{S}_{\lambda}$ is the invere image of the
zero-section.
 (When
$\lambda$
is in a face
$F_P$
of
$\overline{C}$,
we define analogous polynomials
$q_{\lambda}^P(x,t)$
and
$q_{\lambda,b}^P(x)$
by taking the product in (\ref{rep poly}) to be over
$w \in W_P \backslash W.$
These give the reduced equations in this case, and
$q_{\lambda}$
is an appropriate power.)

Over $B_s$ there is a universal $K$-valued cameral cover
$$   \widetilde{\cal S}  \longrightarrow  B_s
$$
with ramification divisor $R \subset \widetilde{\cal S}$. From the relative
Picard,
$$   \mbox{Pic}( \widetilde{\cal S}  /  B_s)
$$
we concoct the relative $N$-shifted, $R$-twisted  Prym
$$   \mbox{Prym}_{\Lambda ,R}( \widetilde{\cal S}  /  B_s).
$$
By Theorem \ref{main}, this can then be considered as a parameter space
$ {\cal M}_{S,G,K} $
for all $K$-valued $G$-principal Higgs bundles on $S$. (Recall that our objects
are assumed to be everywhere {\em regular}!)  It comes with a `Hitchin map',
namely the projection to $B_s$, and the fibers corresponding to smooth
projective $\widetilde{S}$ are abelian varieties. When $S$ is a smooth,
projective curve, we recover this way the algebraic complete integrability of
Hitchin's system and its generalizations.

\section {Symplectic and Poisson structures}
\label{symplectic}\
\indent The total space of  Hitchin's original system  is a cotangent bundle,
hence has a natural symplectic structure. For the polynomial matrix systems of
\cite{B} and \cite{AHH} there is a natural Poisson structure which  one writes
down explicitly.

In \cite{Bn} and \cite{M1}, this result is extended to the systems
${\cal M}_{C,K}$
of $K$-valued GL(n) Higgs bundles on $C$, when
$K \approx \omega_C(D)$
for an effective divisor $D$ on $C$. There is a general-nonsense pairing on the
cotangent spaces, so the point is to check that this pairing is `closed', i.e.
satisfies the identity required for a Poisson structure. Bottacin does this by
an explicit computation along the lines of \cite{B}. Markman's idea is to
consider the moduli space
${\cal M}_D$
of stable vector bundles on $C$ with level-$D$ structure. He realizes an open
subset
${\cal M}^0_{C,K}$
of
${\cal M}_{C,K}$,
parametrizing Higgs bundles whose covers are nice, as a quotient (by an action
of the level group) of
$T^*{\cal M}_D$,
so the natural symplectic form on
$T^*{\cal M}_D$
descends to a Poisson structure on
${\cal M}^0_{C,K}$. This is identified with the general-nonsense form (wherever
both exist),
proving its closedness.

In \cite{Muk}, Mukai constructs a symplectic structure on the moduli space of
simple sheaves on a $K3$ surface $S$. Given a curve
$C \subset S$,
one can consider the moduli of sheaves having the numerical invariants of a
line bundle on a curve in the linear system
$  |nC|  $
on $S$. This has a support map to the projective space
$  |nC|  $,which turns it into an ACIHS. This system specializes, by a
`degeneration to the normal cone' argument, cf. \cite{DEL},  to Hitchin's,
allowing translation of various results about Hitchin's system (such as
Laumon's description of the nilpotent cone, cf. \cite{L} ) to Mukai's.

In higher dimensions, the moduli space $\cal M$ of $\Omega^1$-valued Higgs
bundles carries a natural  symplectic structure \cite{S}. (Corlette points out
in \cite{C} that certain components of  an open subet in $\cal M$ can be
described as cotangent bundles.) It is not clear at the moment exactly when one
should expect to have an ACIHS,  with symplectic, Poisson or quasi symplectic
structure, on the moduli spaces of $K$-valued Higgs bundles for higher
dimensional $S$, arbitrary $G$, and arbitrary vector bundle $K$. A beautiful
new idea \cite{M2} is that Mukai's results extend to the moduli of those
sheaves on a (symplectic, Poisson or quasi symplectic) variety $X$ whose
support in $X$ is {\em Lagrangian.}
Again, there is a general-nonsense pairing. At points where the support is
non-singular projective, this can be identified with another, more geometric
pairing, constructed  using the {\em cubic condition} of \cite{DM1}, which is
known to satisfy the closedness requirement.  This approach is quite powerful,
as it includes many non-linear examples such as Mukai's, in addition to the
line-bundle valued spectral systems of \cite{Bn,M1} and also Simpson's
$\Omega^1$-valued GL(n)-Higgs bundles: just take $X := T^*S
\stackrel{\pi}{\rightarrow} S$, with its natural symplectic form, and the
support in $X$ to be proper over $S$ of degree n; such sheaves correspond to
Higgs bundles by $\pi_*$.

The structure group $GL(n)$ can of course be replaced by an arbitrary reductive
group $G$. Using Theorem \ref{main}, this yields (in the analogous cases) a
Poisson structure on the Higgs moduli space  ${\cal M}_{S,G,K}$  described at
the end of the previous section. The fibers of the generalized Hitchin map are
Lagrangian with respect to this structure. Along the lines of our general
approach, the necessary modifications are clear: $\pi_*$ is replaced by the
equivalence of Theorem \ref{main}. One thus considers only Lagrangian supports
which retain a $W$-action, and only {\em equivariant} sheaves on them (with the
numerical invariants of a line bundle). These two restrictions are symplecticly
dual, so the moduli space of Lagrangian sheaves with these invariance
properties is a symplectic (respectvely, Poisson) subspace of the total  moduli
space, and the fibers of the Hitchin map are Lagrangian as expected.

A more detailed review of the ACIHS aspects of Higgs bundles will appear in
\cite{DM2}.

\section {Some applications and problems}
\label{apps}

\noindent \underline{\bf Some applications}  \nopagebreak

\noindent  In \cite{H1}, Hitchin used his integrable system to compute several
cohomology groups of  the moduli space ${\cal SM}$  (of rank 2, fixed odd
determinant vector bundles on a curve $C$) with coefficients in symmetric
powers of its tangent sheaf ${\cal T}$. The point is that the symmetric algebra
$Sym{\bf \dot{\ }} {\cal T}$
is the direct image of
$ {\cal O}_{T^*{\cal SM}}$,
and sections of the latter all pull back via the Hitchin map $h$ from functions
on the base $B$, since the fibers of $h$ are open subsets in abelian varieties,
and the missing locus has codimension $\geq 2$. Hitchin's system is used in
\cite{BNR} to compute a couple of "Verlinde numbers" for GL(n), namely the
dimensions
$h^0({\cal M}, \Theta) = 1,   \qquad   h^0({\cal SM}, \Theta) = n^g$.
 These results are now subsumed in the general Verlinde formulas, cf.
\cite{F2}, \cite{BL}, and other  references therein.

A pretty application of spectral covers was obtained by Katzarkov and Pantev
\cite{KP2}. Let $S$ be a smooth, projective, complex variety, and
$\rho : \pi_1(S)\longrightarrow G$
 a Zariski dense representation into a  simple $G$ (over $\bf{C}$). Assume That
the $\Omega^1$-valued Higgs bundle $ ( {\cal V}, \phi) $ associated to $\rho$
by Simpson is (regular and) generically semisimple,  so the cameral cover is
reduced. Among other things, they show that $\rho$ factors through a
representation of an orbicurve if and only if  the non-standard component
$Prym_{\epsilon}(\widetilde{S})$
 is non zero, where
$\epsilon $
is the one-dimensional sign representation of $W$.
(In a sense, this is the opposite  of
$Prym_{\Lambda}(\widetilde{S})$:
while
$Prym_{\Lambda}(\widetilde{S})$
is common to
$\mbox{Pic}(\widetilde{S}_P)$
for all proper Weyl subgroups,
$Prym_{\epsilon}(\widetilde{S})$
occurs in none except for the full cameral Picard.)

Another application is in \cite{KoP}: the moduli spaces of SL(n)- or
GL(n)-stable bundles on a curve have certain obvious automorphisms, coming from
tensoring with line bundles on the curve, from inversion, or from automorphisms
of the curve.  Kouvidakis and Pantev use the dominant direct-image maps from
spectral Picards and Pryms to the moduli spaces to show that there are no
further, unexpected automorphisms. This then leads to a `non-abelian Torelli
theorem', stating that a curve is determined by the isomorphism class of the
moduli space of bundles on it.  \\  \mbox{}\\

\noindent \underline{\bf Compatibility?}  \nopagebreak

\noindent  Hitchin's construction \cite {H2} of the projectively flat
connection on the vector bundle of non-abelian theta functions over the moduli
space of curves does not really use much about spectral covers. Nor do other
constructions of Faltings \cite{F} and Witten et al \cite{APW}. Hitchin's work
suggests that the `right'  approach should be based on comparison of the
non-abelian connection near a curve $C$ with the abelian connection for
standard theta functions on spectral covers $\widetilde{C}$ of $C$. One
conjecture concerning the possible relationship between these connections
appears in \cite{A}, and some related versions have been attempted by several
people, so far in vain. What's missing is a compatibility statement between the
actions of the two connections on pulled-back sections.  If the expected
compatibility turns out to hold, it would give another proof of the projective
flatness. It should also imply projective finiteness and projective unitarity
of mo!
nodromy for the non-abelian thetas
, and may or may not bring us closer to a `finite-dimensional' proof of
Faltings' theorem (=the former Verlinde conjecture).\\  \mbox{}\\

\noindent \underline{\bf {Irregulars?} }   \nopagebreak

\noindent The Higgs bundles we consider in this survey are assumed to be
everywhere regular. This is a reasonable assumption for line-bundle valued
Higgs bundles on a curve or surface, but {\em not} in $\dim \geq 3$. This is
because the complement of ${\frak g}_{{reg}}$ has codimension 3 in ${\frak g}$.
The source of the difficulty is that the analogue of (\ref{commutes}) fails
over
${\frak g}$.  There are two candidates for the universal cameral cover:
$\widetilde{\frak g}$, defined by the left hand side of (\ref{commutes}), is
finite over ${\frak g}$ with $W$ action, but does not have a family of line
bundles parametrized by $\Lambda$.
These live on   $\stackrel{\approx}{\frak g}$, the object defined by the right
hand side, which parametrizes pairs
$(x,{\frak b}),    \qquad          x \in {\frak b} \subset {\frak g}$ .
This suggests that the right way to analyze irregular Higgs bundles may involve
spectral data consisting of a tower
$$  \stackrel{\approx}{S}  \stackrel{\sigma}{\longrightarrow}  \widetilde{S}
\longrightarrow  S
$$
together with a homomorphism
$ {\cal L} : \Lambda \longrightarrow  \mbox{Pic}(\stackrel{\approx}{S})$
such that  the collection of sheaves
$$  \sigma_*({\cal L}(\lambda)),   \qquad         \lambda \in \Lambda
$$
 on
$\widetilde{S}$
is
$R$-twisted $W$-equivariant in an appropriate sense.  As a first step, one may
wish to understand the direct images
$ R^i \sigma_*({\cal L}(\lambda))  $
and in particular the cohomologies
$H^i(F, {\cal L}(\lambda))$
where $F$, usually called a {\em Springer fiber}, is a fiber of $\sigma$. For
regular $x$, this fiber is a single point. For $x=0$, the fiber is all of
$G/B$, so the fiber cohomology is given by the Borel-Weil-Bott theorem. The
question may thus be considered as a desired extension of BWB to general
Springer fibers.

\noindent University of Pennsylvnia, Philadelphia, PA 19104-6395\\
Donagi@math.upenn.edu

\end{document}